\documentclass[a4paper,reqno,10pt,oneside]{amsart}
\usepackage{amsmath}
\usepackage[left=2.61cm,right=2.61cm,top=2.72cm,bottom=2.72cm]{geometry}
\usepackage[colorlinks=true,urlcolor=blue,
citecolor=red,linkcolor=blue,linktocpage,pdfpagelabels,
bookmarksnumbered,bookmarksopen]{hyperref}
\usepackage[english]{babel}
\usepackage{graphicx}
\usepackage[toc,page,header]{appendix}
\usepackage{mathrsfs}
\usepackage{amssymb,amsmath, amsfonts, amsthm}
\usepackage{latexsym}
\usepackage{enumitem}
\allowdisplaybreaks[1]
\numberwithin{equation}{section}
\usepackage[numbers]{natbib}

\renewcommand{\(}{\left(}
\renewcommand{\)}{\right)}
\renewcommand{\[}{\left[}
\renewcommand{\]}{\right]}
\newtheorem{theorem}{Theorem}[section]
\newtheorem{proposition}[theorem]{Proposition}
\newtheorem{corollary}[theorem]{Corollary}
\newtheorem{lemma}[theorem]{Lemma}
\theoremstyle{definition}
\newtheorem{remark}[theorem]{Remark}
\theoremstyle{definition}
\newtheorem{definition}[theorem]{Definition}
\theoremstyle{definition}

\renewcommand{\le}{\leqslant}
\renewcommand{\ge}{\geqslant}

\newcommand{\Cl}{{\mathcal L}}

\newcommand{\N}{\mathbb{N}}
\renewcommand{\S}{\mathbb{S}}

\newcommand{\beq}{\begin{equation}}
\newcommand{\eeq}{\end{equation}}
\newcommand{\beqs}{\begin{equation*}}
\newcommand{\eeqs}{\end{equation*}}
\newcommand{\beqn}{\begin{eqnarray}}
\newcommand{\eeqn}{\end{eqnarray}}
\newcommand{\beqns}{\begin{eqnarray*}}
\newcommand{\eeqns}{\end{eqnarray*}}
\newcommand{\bdoc}{\begin{document}}
\newcommand{\edoc}{\end{document}}
\newcommand{\be}{\begin{enumerate}}
\newcommand{\ee}{\end{enumerate}}
\newcommand{\bdescr}{\begin{description}}
\newcommand{\edescr}{\end{description}}
\newcommand{\ba}{\begin{array}}
\newcommand{\ea}{\end{array}}
\newcommand{\intR}{\int_{\mathbb R^N}}
\newcommand{\R}{\mathbb R}
\newcommand{\B}{\mathbb B}
\newcommand{\C}{\mathbb C}
\renewcommand{\H}{\mathcal H}
\renewcommand{\L}{\mathbb L}
\newcommand{\parallelsum}{\mathbin{\!/\mkern-5mu/\!}}
\newcommand{\e}{\varepsilon}
\newcommand{\SD}{\Sigma_D}
 \renewcommand{\(}{\left(}
\renewcommand{\)}{\right)}
\renewcommand{\[}{\left[}
\renewcommand{\]}{\right]}
\renewcommand{\appendixpagename}{\centering Appendix}

\newcommand{\todo}[1]{\text{\colorbox{yellow}{#1}}}

\newcommand{\modified}[1]{\text{\colorbox{green}{#1}}}
\begin{document}
\title[Stability for the Sobolev inequality in cones]{Stability for the Sobolev inequality in cones}

\author{Giulio Ciraolo}
\address{G. Ciraolo. Dipartimento di Matematica ``Federigo Enriques'',
Universit\`a degli Studi di Milano, Via Cesare Saldini 50, 20133 Milano, Italy}
\email{giulio.ciraolo@unimi.it}

\author{Filomena Pacella}
\address{F. Pacella. Dipartimento di Matematica ``Guido Castelnuovo'', Sapienza Universit\`a di Roma,  P.le Aldo Moro 2, 00185 Roma, Italy}
\email{filomena.pacella@uniroma1.it}

\author{Camilla Chiara Polvara}
\address{C.C. Polvara. Dipartimento di Matematica ``Guido Castelnuovo'', Sapienza Universit\`a di Roma,  P.le Aldo Moro 2, 00185 Roma, Italy}
\email{camilla.polvara@uniroma1.it}

\subjclass[49J40, 26D10, 35A23, 35J91, 35B33]{49J40, 26D10, 35A23, 35J91, 35B33}

\begin{abstract}
 We prove a quantitative Sobolev inequality in cones of Bianchi-Egnell type, which implies a stability property. Our result holds for any cone as long as the minimizers of the Sobolev quotient are nondegenerate, which is the case of most cones. When the minimizers are the classical bubbles we have more precise results. Finally, we show that local estimates are not enough to get the optimal constant for the quantitative Sobolev inequality.
\end{abstract}

\maketitle

\section{Introduction}
In this paper, we discuss a Sobolev inequality in cones with the aim of providing a quantitative version of it. 
Let ${D}$ be a smooth domain on the unit sphere $\S^{N-1}$ in $\R^N$, $N\ge 3$ and let $\Sigma_D$ be the cone spanned by $D$, i.e.:
\begin{equation} \label{Sigma_def}
 {\Sigma_D}=\{x\in\R^N; x=sq, s\in (0,+\infty),q\in D\}\,,
\end{equation}
and assume that it is a Lipschitz cone.
We consider the Sobolev space:
$$\mathcal{D}^{1,2}(\Sigma_D)=\{u\in L^{2^*}(\Sigma_D):|\nabla u|\in L^2(\Sigma_D)\},$$
where $2^*=\frac{2N}{N-2}$ is the critical exponent, and the corresponding Sobolev quotient:
\beq\label{Sobquot}
Q_D(u)=\frac{\bigg(\int_{\Sigma_D} |\nabla u|^2\, dx\bigg)^\frac{1}{2}}{\bigg(\int_{\Sigma_D} |u|^{2^*}\, dx\bigg)^\frac{1}{2^*}}, \quad u\in \mathcal{D}^{1,2}(\Sigma_D), u\ne 0. 
\eeq

In \cite[Theorem 3.3 and Proposition 2.7]{CP} it is proved that if $D$ is strictly contained in the hemisphere $\S_+^{N-1}=\S^{N-1}\cap\{x=(x_1,...,x_N,x_N>0)\}$ then the infimum:
\beq\label{best cost}
S_D=\inf_{u\in\mathcal{D}^{1,2}(\Sigma_D)\setminus\{0\}}Q_D(u)
\eeq
is achieved.

Thus from now on we assume that $\bar D\subset \S_+^{N-1}$, so that the following Sobolev inequality holds:
\begin{equation} \label{Sob_ineq}
\|\nabla\varphi\|_2\ge  S_D\|\varphi\|_{2^*}, \quad \forall\varphi\in \mathcal{D}^{1,2}(\Sigma_D),
\end{equation}
where $\|\cdot\|_q$ denotes the norm $\|\cdot\|_{L^q(\Sigma_D)}$ and $S_D$ is the ``best constant'' for it.

When $D=\S^{N-1}$, i.e. the cone is the whole $\R^N$, it is well known that the best constant for the corresponding Sobolev inequality is achieved by the radial functions (usually called bubbles):
\beq\label{bubbles}
U(c,x)=c(1+|x|^2)^{-\frac{N-2}{2}}\quad c\in\R
\eeq
as well as by any translation or rescaling of them. Actually, these functions are the only minimizers of the Sobolev quotient in $\R^N$ and they are also the only positive solutions in $\mathcal{D}^{1,2}(\R^N)$ of the critical equation:
\beq\label{criticaleq}
-\Delta u=u^{2^*-1}\quad \text{ in }\R^N.
\eeq
Since the bubbles are radial functions, they belong to the space $\mathcal{D}^{1,2}(\Sigma_D)$, for any $D\subset\S^{N-1}$, and hence it is natural to ask whether or not they are minimizers for \eqref{best cost}. Moreover, it is well known that they solve the Neumann problem
\beq\label{Neumanpb}
\begin{cases}
    -\Delta u=u^{2^*-1}&   \text{in }\Sigma_D\\
    \frac{\partial u}{\partial\nu}=0 & \text{on }\partial\Sigma_D
\end{cases}
\eeq
and are the only radial solutions of \eqref{Neumanpb}. Thus we wonder if they are the only positive solutions of \eqref{Neumanpb}. The answers to the above questions depend on the cone, i.e. on the domain $D\subset\S^{N-1}_+$ which spans $\Sigma_D$.

In fact, if the cone is convex it can be proved that the only minimizers for \eqref{best cost} are the bubbles. This is a consequence of an isoperimetric inequality obtained in \cite{LP}(see also \cite{RitRos}, \cite{FigIndrei},\cite{CabRosOtSer}) and of the corresponding symmetrization technique in cones (\cite{LPT}, \cite{PT}). The same symmetry result holds in almost convex cones by the extension of the isoperimetric inequality of \cite{BaerFigalli}.

Concerning \eqref{Neumanpb} it was proved in \cite{LPT} that if the cone is convex then the only positive solutions are the bubbles. This result was recently generalized in \cite{CFR}. 

On the contrary break of symmetry results have been proved in \cite{CPP} showing the existence of nonradial minimizer for \eqref{best cost}, and hence also of nonradial positive solutions of \eqref{Neumanpb}, whenever the cone $\Sigma_D$ is such that $\lambda_1(D)<N-1$, where $\lambda_1(D)$ is the first nontrivial Neumann eigenvalue on $D$ of the Laplace-Beltrami operator $-\Delta_{\S^{N-1}}$. Note that for any convex domain $D\subset\S^{N-1}$ it holds $\lambda_1(D)>N-1$.

Hence it is natural to ask whether the condition $\lambda_1(D)\ge N-1$ would represent a threshold to get radial symmetry for \eqref{best cost} and \eqref{Neumanpb}. This is an open question and we conjecture that the answer should be affirmative.

Thus for the Sobolev inequality \eqref{best cost} we have two different scenarios: either the only minimizers are bubbles or there are nonradial minimizers. In the last case, there could be even infinitely many nonradial ``independent'' minimizers, where by the word independent we mean not obtained by one fixed minimizer by rescaling or multiplication by a constant. Let us observe that since the cone $\Sigma
_D$ has its vertex at the origin and $D$ is smooth we do not have invariance by translation. If the cone was not smooth, for example, if $\Sigma_D=\R^{N-k}\times C$ where $C\subset\R^k$ is a cone which does not contain straight lines, then there would be invariance with respect to suitable translations. This does not add particular difficulties to the proof of our results and we prefer to neglect this case. 

The aim of the present paper is to show a quantitative version of the Sobolev inequality \eqref{Sobquot}. Inspired by a famous inequality obtained by Bianchi and Egnell in \cite{BE} in $\R^N$ we aim at estimating the Sobolev ``deficit'' for a generic function $\varphi\in\mathcal{D}^{1,2}(\Sigma_D)$ in terms of the distance of $\varphi$ from the minimizers. 

To be more precise let us consider a fixed minimizer $V$ for \eqref{Sob_ineq} such that $\|\nabla V\|_2=1$ and the rescaled functions:
\beq\label{Vminimizer}
V_s(x)=sV(s^\frac{2}{N-2}x), \quad s\in (0,+\infty).
\eeq
Obviously $V_s$ are still minimizers for \eqref{Sob_ineq} and so it is any function $cV_s, c\in\R\setminus\{0\}$. Thus we define the 2-dimensional manifold:
$$
\mathcal{M}_V=\{cV_s(x), c\in\R\setminus\{0\}, s\in\R^+\}
$$
and the distance of a function $\varphi\in\mathcal{D}^{1,2}(\Sigma_D)$ from $\mathcal{M}_V$:
$$
d(\varphi,\mathcal{M}_V)=\inf_{u\in \mathcal {M_V}}\|\nabla(\varphi-u)\|_{L^2({\Sigma_D})}=\inf_{c\in\R\setminus\{0\}, s\in\R^+}\|\nabla(\varphi-cV_s)\|_{L^2({\Sigma_D})}.
$$
In what follows we consider nondegenerate minimizers whose definition is given in Section 2, (Definition \ref{defnondegen}). Our main result is the following
\begin{theorem} \label{thm_main}
Let ${\Sigma_D}$ be a cone as in \eqref{Sigma_def}, with $\bar D\subset\S^{N-1}_+$ and assume that every minimizer for \eqref{best cost} is nondegenerate.
Then 
\begin{enumerate}
    \item up to rescaling and multiplication by a constant, there exists only a finite number of minimizers that we denote by $V^{(1)},...,V^{(m)}, m\in\N^+$ with $\|\nabla V^{(i)}\|_2=1, i=1,...,m$.
    \item there exists a positive constant $C_D$ such that the following  quantitative Sobolev inequality holds
\beq\label{quant_sobolev0}
\|\nabla\varphi\|_2^2- S_D^2\|\varphi\|_{2^*}^2\ge C_D\bigg (\min_{i=1,...,m}d(\varphi,\mathcal{M}_{V^{(i)}})^2\bigg), 
\eeq
for every $\varphi\in\mathcal{D}^{1,2}(\Sigma_D)\setminus\{0\}$. 
\end{enumerate}
\end{theorem}
It is important to stress that the above inequality is valid whatever the minimizers are, as long as they are nondegenerate, which is the case of most cones, since generically all eigenvalues of the associated operator $\mathcal{L}_V$ are simple (see \cite{Henry}). In the case when the minimizers are the bubbles then their nondegeneracy holds in any cone $\Sigma_D$ for which $\lambda_1(D)>N-1$, see Corollary \ref{corollary}. The previous theorem can be considered as a stability result for the Sobolev inequality \eqref{Sob_ineq}, because it claims that if the ``deficit'' $\|\nabla\varphi\|_2^2-\hat S_D^2\|\varphi
\|_{2^*}^2$ is very small for some $\varphi\in\mathcal{D}^{1,2}(\Sigma_D)$ then $\varphi$ must be very close to one of the minimizers and the distance from it can be estimated uniformly by the deficit.

The proof of Theorem \ref{thm_main} follows the approach of \cite{BE}. To quantify the Sobolev deficit we first get an expansion of it for functions which are close to some of the manifolds of minimizers. This works for nondegenerate local minimizers, using also the spectral analysis of the operators $\mathcal{L}_{V^{(i)}}$. This provides a local estimate of the Sobolev deficit involving the eigenvalues of $\mathcal{L}_{V^{(i)}}$. Once we have this the conclusion follows by a compactness argument.

The local estimate also provides an upper bound for the constant $C_D$ in \eqref{quant_sobolev0} which we can assume to be the optimal one, i.e.
\beq\label{bestconstdefic}
C_D=\inf_{\varphi\in \mathcal{D}^{1,2}(\Sigma_D)\setminus \bigcup_{i=1}^m \mathcal{M}_{V^{(i)}}}\frac{\|\nabla \varphi\|_2^2-\hat S_D^2\|\varphi\|_{2^*}^2}{\min_{i=1,...,m}d(\varphi,\mathcal{M}_{V^{(i)}})^2}
\eeq
Indeed by Lemma \ref{lemma1} we get:
\beq\label{Cd}
C_D\le\min_{i=2,...,m}\bigg(1-\frac{\mu_2^i}{\mu_3^i}\bigg)\,.
\eeq
where $\mu_k^i, k\in\N$ are the eigenvalues of the weighted operators $\mathcal{L}_{V^{(i)}}$. In the case when the only minimizers for \eqref{best cost} are the bubbles  then \eqref{Cd} becomes:
\beq\label{cd*}
C_D\le C_D^*
\eeq
where $C_D^*$ is an explicit constant computed in Corollary \ref{corineq}.

When $D=\S^{N-1}$ then $C_D^*=\frac{4}{N+4}$ and it is proved in \cite{Konig1} that the inequality \eqref{cd*} is strict.
In Section 3 we get a similar result showing that \eqref{Cd} is strict when  $N=3,4,5$ or when $N\ge 6$, the only minimizers for \eqref{best cost} are the bubbles, $\lambda_1(D)\ge2N$ and $D$ has some symmetries (see Theorem \ref{strict}).

Moreover, again for $D=\S^{N-1}$, it is proved in \cite{Konig2} that $C_D$ is achieved. We believe that this is true also in any cone $\Sigma_D$, whenever the minimizers are nondegenerate and we plan to study this question in a future research.

Finally, interesting estimates for $C_D$ have been proved in \cite{Estebanet} when $D=\S^{N-1}$, obtaining also an alternative proof of the Bianchi Egnell result. However, these estimates strongly exploit the explicit form of the bubbles and their radial symmetry. Some similar result should be obtained when the minimizers for \eqref{best cost} in $\Sigma_D$ are the bubbles, though the proof could be much more involved since $D$ is not $S^{N-1}$ but any domain on $S^{N-1}_+$.
\medskip

\noindent {\bf  Organization of the paper.} In Section \ref{spectan} we study the spectrum of the weighted operator $\mathcal{L}_V$ and prove some preliminary results both for a generic minimizer $V$ and for the case when the minimizers are the bubbles. In Section \ref{proofth1} we prove a local estimate of the Sobolev deficit from which the proof of Theorem \ref{thm_main} derives. In the same Section, we show that the strict inequality in \eqref{cd*} holds for some symmetric cones. Finally, a short Appendix explains how to compute the radial eigenvalues of  $\mathcal{L}_V$, when $V$ is a bubble.
\medskip

\noindent{\bf Acknowledgements.} The authors thank Bruno Premoselli for interesting conversations which helped to improve our results.
Research partially supported by Gruppo Nazionale per l'Analisi Matematica, la Pro\-ba\-bi\-li\-t\`a e le loro Applicazioni (GNAMPA) of the Istituto Nazionale di Alta Matematica (INdAM). G. C. has been supported by the Research Project of the Italian Ministry of University and Research (MUR) Prin 2022 “Partial differential equations and related geometric-functional inequalities”, grant number 20229M52AS004. F.P. and C.C.P. have been supported by PRIN 2022 PNRR project 2022AKNSE4-Next Generation EU “Variational and Analytical aspects of Geometric PDEs”, founded by the European Union.

\section{Spectral analysis}\label{spectan}

\subsection{The case of a generic local minimizer}

Let $V$ be a local minimizer of the Sobolev quotient \eqref{Sobquot} and consider the space $L^2(V^{2^*-2})$, that is given by the functions $v$ such that 
$$
\int_{{\Sigma_D}} V^{2^*-2} v^2<\infty.
$$
From \cite[Proposition 3.1]{CPP}, we know that $V$ is bounded, of class $C^2(\Sigma_D\setminus\{0\})$ and  $V\le C |x|^{2-N}$ as $|x|\to\infty$.

We start with a preliminary result.

 \begin{proposition}[Hardy inequality] \label{prop_Hardybe}
Let $N \geq 3$  and let ${\Sigma_D}$ be given by \eqref{Sigma_def}. For any $u\in H_0^1({\Sigma_D}\cup \Gamma_1 )$, and for any $\epsilon\ge 0$  it holds
\begin{equation} \label{hardybe}
\frac{(N-4+2\epsilon)^2}{4}\int_{\Sigma_D}\frac{u^2}{|x|^{4-2\epsilon}}\, dx\le \int_{\Sigma_D}\frac{|\nabla u|^2}{|x|^{2-2\epsilon}} \, dx,
\end{equation}
\end{proposition}
\begin{proof}
 We apply a classical argument used to prove Hardy-type inequalities (see for example \cite{CW}). 
Given $u\in C_c^\infty({\bar\Sigma_D}\setminus\{0\})$, for any $\epsilon\ge 0$ it holds
$$
\int_{\Sigma_D} \bigg|\frac{\nabla u}{|x|^{1-\epsilon}}+t\frac{x}{|x|^{3-\epsilon}}u\bigg|^2\, dx\ge 0
$$
that is
\beq\label{hardystep1}
0\le \int_{\Sigma_D}\frac{|\nabla u|^2}{|x|^{2-2\epsilon}} \, dx+t^2 \int_{\Sigma_D} \frac{u^2}{|x|^{4-2\epsilon}}\, dx+2t\int_{\Sigma_D} u\frac{x\cdot\nabla u}{|x|^{4-2\epsilon}}\, dx.
\eeq
Now we consider 
\begin{multline}\label{Hardystep2}
2\int_{\Sigma_D} u\frac{x\cdot\nabla u}{|x|^{4-2\epsilon}}\, dx=\int_{\Sigma_D} \text{div}  \bigg(u^2\frac{x}{|x|^{4-2\epsilon}}\bigg)\, dx-\int_{\Sigma_D} u^2\text{div} \bigg(\frac{x}{|x|^{4-2\epsilon}}\bigg)\, dx\\
=0-(N-4+2\epsilon)\int_{\Sigma_D}\frac{u^2}{|x|^{4-2\epsilon}}\, dx.
\end{multline}
Substituting \eqref{Hardystep2} in \eqref{hardystep1} we get
$$
0\le \int_{\Sigma_D}\frac{|\nabla u|^2}{|x|^{2-2\epsilon}} \, dx+t^2 \int_{\Sigma_D} \frac{u^2}{|x|^{4-2\epsilon}}\, dx-(N-4+2\epsilon)t\int_{\Sigma_D}\frac{u^2}{|x|^{4-2\epsilon}}\, dx,
$$
that is
$$
(N-4+2\epsilon-t)t\int_{\Sigma_D}\frac{u^2}{|x|^{4-2\epsilon}}\, dx\le \int_{\Sigma_D}\frac{|\nabla u|^2}{|x|^{2-2\epsilon}} \, dx\,.
$$
An optimization it $t$ yields that the best possible choice of the constant on the l.h.s. is given by 
$$
\frac{(N-4+2\epsilon)^2}{4}\int_{\Sigma_D}\frac{u^2}{|x|^{4-2\epsilon}}\, dx\le \int_{\Sigma_D}\frac{|\nabla u|^2}{|x|^{2-2\epsilon} }\, dx,
$$
and we conclude the proof.
\end{proof}
Next we study the spectrum of the operator $\Cl_{V}=V^{2-2^*}\Delta$ on $L^2(V^{2^*-2})$.

\begin{lemma}\label{lemma_spettro}
$\Cl_{V}=V^{2-2^*}\Delta$ on $L^2(V^{2^*-2})$ has discrete spectrum.
\end{lemma}

\begin{proof}
This follows from the fact that the embedding $\mathcal{D}^{1,2}(\Sigma_D)\hookrightarrow L^2(V^{2^*-2})$ is compact, and hence the spectrum is discrete, which can be proved as in \cite[Corollaries 6.2 and 6.3]{FN}. We notice that in \cite[Corollaries 6.2]{FN}, Rellich-Kondrachov compact embedding theorem and the Hardy type inequality are used. In our setting, the Rellich-Kondrachov Theorem holds, while the validity of Hardy inequality is proved in Proposition \ref{prop_Hardybe}.  For $N\ne 4$ we take $\epsilon=0$ in Proposition \ref{prop_Hardybe} and proceed in the same way as in \cite[Corollaries 6.2]{FN} noticing that $V\in L^1_{loc}(\Sigma_D)$ and $V\le C |x|^{2-N}$ for $|x|\to\infty$ with $C>0$ (see \cite{CP}). For $N=4$ we take $\epsilon>0$ and we have:
\begin{multline*}
\int_{{\Sigma_D}\setminus B_k}\frac{u^2}{|x|^4}\le \frac{1}{k^{2\epsilon}}\int_{{\Sigma_D}\setminus B_k}\frac{u^2}{|x|^{4-2\epsilon}}\le \frac{4}{(N-4+2\epsilon)^2} \frac{1}{k^{2\epsilon}}\int_{{\Sigma_D}\setminus B_k}\frac{|\nabla u|^2}{|x|^{2+\epsilon}} \, dx\\
\le \frac{1}{\epsilon^2k^{2\epsilon}} \int_{{\Sigma_D}\setminus B_k}\frac{|\nabla u|^2}{|x|^{2-2\epsilon}} \, dx\le  \frac{1}{\epsilon^2k^{2}} \int_{{\Sigma_D}\setminus B_k}|\nabla u|^2 \, dx\, .
\end{multline*}
Then we conclude as in \cite[Corollaries 6.2 and 6.3]{FN}.
\end{proof}
Let $\N_+$ be the set of positive natural numbers, we denote by $\mu_i$, $i\in\N_+$ the eigenvalues of $\Cl_{V}$, i.e.
\begin{equation} \label{pb_mu_def}
\begin{cases}
-V^{2-2^*}\Delta v = \mu_i v & \text{in } {\Sigma_D} \\
\partial_\nu v = 0 & \text{on } \partial {\Sigma_D} \\
v \in \mathcal{D}^{1,2}(\Sigma_D). &
\end{cases}
\end{equation}
We notice that the eigenvalues of \eqref{pb_mu_def} are the same if, instead of considering $V$, we consider a rescaled minimizer $V_s$. Indeed
let $\mu_s$ be an eigenvalue of \eqref{pb_mu_def} corresponding to $V_s$ with $s$ fixed and $\psi(x)$ being the related eigenfunction, namely
$$
-\Delta \psi = \mu_s V_s^{2^*-2} \psi.
$$
We take as a coordinate $y=s^{-\frac{2}{N-2}}x$, then we call $\tilde\psi(x):=\psi(y)$ and we see that $\tilde\psi$ satisfies
$$
-\Delta \tilde\psi(x) = -s^{-\frac{4}{N-2}}\Delta \psi(s^{-\frac{2}{N-2}}x)= \mu_s V(s^{\frac{2}{N-2}}s^{-\frac{2}{N-2}}x)\psi(s^{-\frac{2}{N-2}}x)=\mu_sV^{2^*-2}(x)\tilde \psi(x).
$$
Then we can conclude that the eigenvalues do not depend on $s$.

Thus we can assume that $s=1$ and let $V_{s}=V$. 

\begin{proposition}\label{proposspettro}
Let $V$ be a local minimizer of the Sobolev quotient, and let 
$$
S_V:=\frac{\|\nabla V\|_{L^2({\Sigma_D})}}{\|V\|_{L^{2^*}({\Sigma_D})}}.
$$
Then for the eigenvalue problem \eqref{pb_mu_def}, we have 
\begin{itemize}
\item[(i)] $\mu_1=S_V^{2^*}$ is simple with eigenfunction $V$; 
\item[(ii)] $\mu_2=(2^*-1)S_V^{2^*}$ and a corresponding eigenfunction is $\partial_s V_{s}|_{s=1}$ ($V_s$ as in \eqref{Vminimizer}). 
\end{itemize}
\end{proposition}

\begin{proof}
We can suppose w.l.o.g. that $\|\nabla V\|_{\Sigma_D}^2=1$.
    We compute the second variation of
  $$
\mathcal{F}(u) = \frac{\|\nabla u\|_{L^2(\Sigma_D)}}{\|u\|_{L^{2^*}(\Sigma_D)}},
$$
in $V$. We have that for any $\eta\in \mathcal{D}^{1,2}(\Sigma_D)$ and $t\in\R$
$$
\int_{\Sigma_D}(V+t\eta)^{2^*}=\int_{\Sigma_D}(V)^{2^*}+2^*t\int_{\Sigma_D}V^{2^*-1}\eta+2^*(2^*-1)t^2\int_{\Sigma_D}V^{2^*-2}\eta^2+o(t^2),
$$
and thus
\begin{multline*}
    \left( \int_{\Sigma_D} (V + t\eta)^{2^*} \right)^{-\frac{2}{2^*}}
= \left( \int_{\Sigma_D} (V)^{2^*} \right)^{-\frac{2}{2^*}}
\Bigg( 
1 - \frac{2t \int_{\Sigma_D} V^{2^*-1} \eta}{\int_{\Sigma_D} (V)^{2^*}}
- \frac{t^2 (2^*-1) \int_{\Sigma_D} V^{2^*-2} \eta^2}{\int_{\Sigma_D} (V)^{2^*}}\\
+ \frac{2^* \left(1 + \frac{2}{2^*}\right) t^2 \left(\int_{\Sigma_D} V^{2^*-1} \eta\right)^2}{\left(\int_{\Sigma_D} (V)^{2^*}\right)^2} 
\Bigg)+o(t^2)\,,
\end{multline*}
i.e. 
\begin{multline}\label{comp111}
 \left( \int_{\Sigma_D} (V + t\eta)^{2^*} \right)^{-\frac{2}{2^*}}
= S_V^2
\Bigg( 
1 -S_V^{2^*}\bigg(2t \int_{\Sigma_D} V^{2^*-1} \eta+ t^2 (2^*-1) \int_{\Sigma_D} V^{2^*-2} \eta^2\bigg)
+\\+ S_V^{2\cdot2^*}2^* \left(1 + \frac{2}{2^*}\right) t^2 \left(\int_{\Sigma_D} V^{2^*-1} \eta\right)^2
\Bigg)+o(t^2).
\end{multline}
Since
\beq\label{comp222}
\|\nabla (V+t\eta)\|_2^2=1+2t\int_{\Sigma_D}\nabla V\cdot\nabla\eta+t^2\int_{\Sigma_D}|\nabla \eta|^2
\eeq
from \eqref{comp111} and \eqref{comp222} we obtain 
\begin{multline*}
\|\nabla (V+t\eta)\|_2^2\|V+t\eta\|_{2^*}^{-\frac{2}{2^*}}=S_V^2+2tS_V^2\bigg(\int_{\Sigma_D}\nabla V\cdot\nabla\eta-S_V^{2^*} \int_{\Sigma_D} V^{2^*-1} \eta \bigg)+t^2S_V^2\bigg\{\int_{\Sigma_D}|\nabla\eta|^2+\\
-4S_V^{2^*} \int_{\Sigma_D} V^{2^*-1} \eta\int_{\Sigma_D}\nabla V\cdot\nabla\eta+S_V^{2^*}\bigg(-(2^*-1)\int_{\Sigma_D} V^{2^*-2} \eta^2+S_V^{2^*}(2^*+2)\left(\int_{\Sigma_D} V^{2^*-1} \eta\right)^2\bigg)\bigg\}+o(t^2).
\end{multline*}
Thus for the first variation $\mathcal{F}'(V)$ we have
\beq\label{firstvar}
\langle \mathcal{F}'(V), \eta\rangle= 2S_V^2\bigg(\int_{\Sigma_D}\nabla V\cdot\nabla\eta-S_V^{2^*} \int_{\Sigma_D} V^{2^*-1} \eta \bigg)=0 \quad \forall\eta \in \mathcal{D}^{1,2}(\Sigma_D)
\eeq
because $V$ is a local minimizer.
The second variation is given by
\begin{multline}\label{second var}
Q_{V}(\eta):=\int_{\Sigma_D}|\nabla\eta|^2-4S_V^{2^*} \int_{\Sigma_D} V^{2^*-1} \eta\int_{\Sigma_D}\nabla V\nabla\eta+
\\
+S_V^{2^*}\bigg(-(2^*-1)\int_{\Sigma_D} V^{2^*-2} \eta^2+S_V^{2^*}(2^*+2)\left(\int_{\Sigma_D} V^{2^*-1} \eta\right)^2\bigg)
\end{multline}
and we have $Q_{V}(\eta)\ge 0$ for any $\eta\in \mathcal{D}^{1,2}(\Sigma_D)$ because $V$ is a local minimizer. From \eqref{firstvar} we get
\beq\label{firstvar2}
\int_{\Sigma_D}\nabla V\cdot\nabla\eta=S_V^{2^*} \int_{\Sigma_D} V^{2^*-1} \eta\,,
\eeq
so, by substituting in \eqref{second var} and by using $Q_{V}(\eta)\ge 0$, we obtain
\beq\label{second var2}\int_{\Sigma_D}|\nabla\eta|^2+(2^*-2)S_V^{2^*} \bigg(\int_{\Sigma_D}\nabla V\cdot\nabla\eta\bigg)^2-S_V^{2^*}(2^*-1)\int_{\Sigma_D} V^{2^*-2} \eta^2\ge 0\quad\forall
\eta\in\mathcal{D}^{1,2}(\Sigma_D)\,.
\eeq
From $\eqref{firstvar2}$ we have that $\mu=S_V^{2^*}$ is an eigenvalue, with eigenfunction $V$. Since $V>0$, from the variational characterization of eigenvalues we deduce that $\mu=S^{2^*}_V$ is the first eigenvalue. Indeed all other eigenfunctions must be orthogonal to $V$ and, therefore, they must be sign changing.
Note that $Q_V(V)=0$, as it follows from \eqref{second var}.
Observe that since each rescaled function $V_s$ satisfies
$$
-\Delta V_s=S_{V_s}^{2^*} V_s^{2^*-1},
$$
by differentiating with respect to $s$ we get
$$
-\Delta \partial_sV_s=S_V^{2^*}(2^*-1) V_s^{2^*-2}\partial_sV_s,
$$
because $S_{V_s}=S_V$.
This implies that $\partial_s V_s$ is an eigenfunction with eigenvalue $S_V^{2^*}(2^*-1)$ and $Q_V(\partial_sV_s)=0$ as it is easy to deduce from \eqref{second var}.
Suppose that there exists an eigenvalue $\bar\mu$ such that $S_V^{2^*}<\bar\mu<(2^*-1)S_V^{2^*}$ then there exists also an eigenfunction $\bar\psi$ orthogonal to the first eigenfunction $V$ such that 
$$
-\Delta\bar\psi=\bar\mu V^{2^*-2}\bar\psi.
$$
From \eqref{second var2} we must have
$$
\int_{\Sigma_D}|\nabla\bar\psi|^2+(2^*-2)S_V^{2^*} \bigg(\int_{\Sigma_D}\nabla V\cdot\nabla\bar\psi\bigg)^2-S_V^{2^*}(2^*-1)\int_{\Sigma_D} V^{2^*-2} \bar\psi^2\ge 0
$$
but, from the orthogonality condition and the definition of $\bar\psi$, this is equivalent to
$$
\bigg(\bar \mu-S_V^{2^*}(2^*-1)\bigg)\int_{\Sigma_D} V^{2^*-2} \bar\psi^2\ge 0\,,
$$
which is true if and only if $\bar \mu\ge (2^*-1)S_V^{2^*}$, a contradiction. Then we conclude that $\mu_2=S_V^{2^*}(2^*-1)$ and a second eigenfunction is given by $\partial_s V_s$.
\end{proof}
\begin{definition}\label{defnondegen}
We say that a local minimizer $V$ is \textbf{nondegenerate} if the second variation $Q_V$ vanishes only on the $2-$dimensional subspace spanned by $\{V,\partial_sV_s|_{s=1}\}$. Equivalently, $Q_V(\psi)>0$ for any $\psi$ orthogonal to $\{V,\partial_sV_s|_{s=1}\}$ in $\mathcal{D}^{1,2}(\bar \Sigma_D)$.
\end{definition}
\begin{remark}
    If $V$ is nondegenerate, then the second eigenvalue $\mu_2$ is simple. Indeed, consider the second eigenfunction $\partial_s V_s|_{s=1}$, then $Q_{V}( \partial_s V_s|_{s=1})=0$ as it is easy to see from \eqref{second var}. Suppose that there exists another independent eigenfunction $\psi$, with eigenvalue $(2^*-1)S_V^{2^*}$, then $\psi$ is orthogonal to the first eigenfunction $V$ and 
    $$Q_{V}( \psi)=\int_{\Sigma_D}|\nabla\psi|^2-(2^*-1)S_V^{2^*}\int_{\Sigma_D} V^{2^*-2} \eta^2=0,$$ 
    hence $V$ would be degenerate.
\end{remark}

\subsection{The case of the bubble}

We consider the case in which the local minimizer is the bubble $U$:
\beq\label{bubbledef}
U(x):=k_0(1+|x|^2)^{-\frac{N-2}{2}}
\eeq
where $k_0$ is a constant such that $\|\nabla U\|_2=1$ (see Remark \ref{k0} below). 
The rescaled bubbles will be denoted by $$
U_s(x):=s U(s^{\frac{2}{N-2}}x).
$$
To describe the spectrum of the operator $\mathcal{L}_U$ we need to know the Neumann eigenvalues $\lambda_j(D), j\in \N$, of the Laplace-Beltrami operator $-\Delta_{\S^{N-1}}$ on $D$, i.e. we consider the following problem:
\begin{equation}\label{pb_lambda_def1}
\begin{cases}
-\Delta_{\S^{N-1}}Y_j = \lambda_j Y_j & \text{on } D\\
\partial_{\nu_D} Y_j = 0 & \text{on } \partial D \,.
\end{cases}
\end{equation}
It is well known that $(-\Delta_{\S^{N-1}})^{-1}$ is compact and selfadjoint in $L^2(D)$ and admits a sequence of eigenvalues 
\begin{equation} \label{lambdabe}
0=\lambda_0<\lambda_1(D)\le \ldots \lambda_j \leq \ldots
\end{equation} 
and corresponding eigenfunctions $Y_j(\theta) \in L^2(D)$ (where $\theta$ is the system of coordinates on $D$ induced by the spherical coordinates in $\R^N$) which form a Hilbert basis for $L^2(D)$ and such that
\beq\label{def:lapleigbe}
-\Delta_{\S^{N-1}}Y_j(\theta)=\lambda_j Y_j(\theta),\quad \theta\in D, 
\eeq
and 
$$
\int_{D} \nabla_\theta Y_j(\theta) \cdot \nabla_\theta Y_i(\theta) d\sigma(\theta) = 0 \quad \text{ for } i \neq j \,,
$$
where we denote by $\nabla_\theta$ the gradient with respect to $\theta \in D$.

In particular, the first eigenvalue is $\lambda_0=0$ and the corresponding eigenfunction is constant.  The second eigenvalue $\lambda_1(D)>0$ has eigenspace generated by $Y_1$.

We denote by $S_U$ the Sobolev quotient $Q_D(U)$ (see \eqref{Sobquot}).

In what follows we consider problem \eqref{pb_mu_def} with $U$ instead of a generic local minimizer $V$
\begin{equation}\label{pb_mu_def2} 
\begin{cases}
\Delta v+\mu U^{2^*-2}v=0 &\text{in } {\Sigma_D}\\
\partial_\nu v=0 &\text{on }  \partial{\Sigma_D}\\
v\in \mathcal{D}^{1,2}(\Sigma_D). &
\end{cases}
\end{equation}
\begin{proposition}\label{separation}
Let $\mu$ be an eigenvalue of \eqref{pb_mu_def2}, and let $\psi$ be a corresponding eigenfunction. Then $\psi=R(r)Y_j(\theta)$ where $Y_j(\theta), j\in\N$ is an eigenfunction of \eqref{pb_lambda_def1} corresponding to $\lambda_j(D)$ and $R$ is a solution to:
\begin{equation} \label{SL_eq1}
\begin{cases}
(r^{N-1}R')'+r^{N-1}(-\lambda_j(D) r^{-2}+ \mu U^{2^*-2})R=0 & \text{in }\R_+ \\
(r^{N-1}R'(r))(0)=0 & \\ 
R(r)=o(r^{2-N+\epsilon})  \text{ as }r\to + \infty \, &
\end{cases}
\end{equation}
for all $\epsilon>0$.
Conversely, if $R$ is a solution to \eqref{SL_eq1} for some $\lambda_j(D)$ and $ \mu$, then $R(r)Y_j(\theta)$ is an eigenfunction of\eqref{pb_mu_def2} with $\mu$ as corresponding eigenvalue.
\end{proposition}

\begin{proof}
Let $\mu$ be fixed and let $ \psi\in \mathcal{D}^{1,2}(\Sigma_D)$ be an eigenfunction associated to $\mu$, i.e. $ \psi $ is a solution to \eqref{pb_mu_def2}. Since $\{Y_j\}$ is an orthonormal base on $L^2(D)$, by using polar coordinates and separation of variables we can write
$$
\psi(r,\theta)=\sum_{j=0}^\infty \hat\psi_j(r)Y_j(\theta)
$$
for $r\in(0,\infty)$ and $\theta\in D\subset \S^{N-1}$, where
\beq\label{eq:eigenfbe}
\hat\psi_j(r):=\int_{D}\psi(r,\theta)Y_j(\theta)d\sigma(\theta) \,.
\eeq

Since $ \psi \not\equiv 0$, there exists $j\in\N$ such that $\hat\psi_j(r)\not \equiv 0$ and we can write
$$
\int_0^\infty r^{N-1}\hat\psi_j'\varphi'dr=\int_0^\infty\int_D r^{N-1} \psi' Y_j(\theta)\varphi' drd\sigma(\theta)=\int_0^\infty\int_D r^{N-1} \psi'\bigg(Y_j(\theta)\varphi\bigg)' drd\sigma(\theta)
$$
for every $\varphi\in \mathcal{D}^{1,2}(\Sigma_D)$. We recall that $\nabla  \psi=(\psi', r^{-2} \nabla_\theta \psi)$ and hence
$$
\nabla\psi\cdot\nabla\bigg(Y_j(\theta)\varphi\bigg)=\psi'\bigg(Y_j(\theta)\varphi\bigg)'+r^2\frac{1}{r^4}\nabla_\theta\psi\cdot \nabla_\theta (Y_j(\theta)\varphi) 
$$
since $ \psi$ is a solution to $\eqref{pb_mu_def2}$ and $U$ is radial, we have
\begin{multline*}
\int_0^\infty r^{N-1}\hat\psi_j'\varphi'dr=-\int_0^\infty\int_D r^{N-3}\nabla_\theta\psi\cdot \nabla_\theta (Y_j(\theta)\varphi) drd\sigma(\theta) \\
+\mu\int_0^\infty\int_D r^{N-1}U^{2^*-2}\psi Y_j(\theta)\varphi drd\sigma(\theta).
\end{multline*}
From \eqref{eq:eigenfbe}, \eqref{pb_lambda_def1} and by using integration by parts, we obtain
\begin{flalign*}
\int_0^\infty r^{N-1}\hat \psi_j'\varphi'dr=&-\int_0^\infty r^{N-3}dr \int_D \nabla_\theta( \psi\varphi)\cdot \nabla_\theta Y_j(\theta) d\sigma(\theta)+ \mu\int_0^\infty r^{N-1}U^{2^*-2}\hat \psi_j\varphi dr\\
=&-\lambda_j(D)\int_0^\infty r^{N-3}\hat \psi_j\varphi dr+ \mu\int_0^\infty r^{N-1}U^{2^*-2}\hat \psi_j\varphi dr.
\end{flalign*}
Hence we have proved that $\hat \psi_j$ is a weak solution to 
\beq
-(r^{N-1}(\hat\psi_j)')'- \mu r^{N-1}U^{2^*-2}\hat\psi_j=-\lambda_j(D) r^{N-3}\hat\psi_j\qquad  \text{ for } r\in(0,\infty).
\eeq
Hence we set $R=\hat \psi_j$ and we need to prove that the boundary conditions in \eqref{SL_eq1} hold. 
To prove that $R(r)=o(r^{2-N})$ as $r\to\infty$ we refer to \cite[Lemma 6.6]{FN}. In this lemma, the authors prove that for any $\epsilon>0$
$$
R(r)\le Cr^{2-N+\epsilon}
$$
for any $r>r_0$ with $r_0>0$ fixed.

We are left to prove that $(r^{N-1}R'(r))(0)=0$ for which we refer to the proof of \cite[Lemma 6.4, (1)]{FN}.

This concludes the proof of the first implication.

The reverse implication holds as well. Namely, if we consider the problem \eqref{SL_eq1} with eigenvalue $ \mu$ with associated eigenfunction $ \psi^{rad}\in \mathcal{D}^{1,2}(\Sigma_D)$ and for some $\lambda_j(D)$, then the function 
$$
\Psi:= \psi^{rad}(r)Y_j(\theta)
$$ 
is such that $\Psi \in \mathcal{D}^{1,2}(\Sigma_D)$.  

Moreover $\Psi$ weakly solves \eqref{pb_mu_def2} corresponding to $\mu$. Indeed for any $\varphi\in \mathcal{D}^{1,2}(\Sigma_D)$ we have
\begin{multline*}
\int_{{\Sigma_D}}\nabla\Psi\cdot\nabla\varphi dx=\int_{{\Sigma_D}}\Psi'\varphi' +\frac{1}{r^2}\nabla_\theta\Psi\cdot\nabla_\theta\varphi dx=\\
=\int_D Y_j(\theta)d\sigma(\theta)\int_0^1r^{N-1}(\psi^{rad})' \varphi'dr+\int_0^1 r^{N-3}\psi^{rad}dr\int_D\nabla_\theta Y_j \cdot\nabla_\theta\varphi d\sigma(\theta)
\end{multline*}
and from \eqref{SL_eq1} and \eqref{def:lapleigbe} it follows
\begin{multline*}
\int_{{\Sigma_D}}\nabla\Psi\cdot\nabla\varphi dx=\int_DY_j(\theta)d\sigma(\theta)\int_0^\infty r^{N-1}( \mu U^{2^*-2}-\lambda_j(D) r^{-2})\psi^{rad}\varphi dr+\\
+\int_0^\infty r^{N-3}\psi^{rad}dr\lambda_j(D)\int_D Y_j(\theta )\varphi d\sigma(\theta)
=\int_{{\Sigma_D}} \mu U^{2^*-2}\Psi\varphi\,dx,
\end{multline*}
and then $\Psi$ weakly solves \eqref{pb_mu_def2} corresponding to $ \mu$.

The fact that $\partial_\nu \Psi=0$ on $\partial{\Sigma_D}$ directly follows from the condition $\partial_\nu Y_j(\theta)=0$ on $\partial D$.
\end{proof}

\begin{proposition} \label{lemma_autov}

Let $D\subset\S^{N-1}$ and $k\in\N_+$ be such that
\beq\label{lemma_autov_cond}
(k-1)(k+N-3)<\lambda_1(D)<k(k+N-2).
\eeq
Then each eigenvalue $\mu_i$ of \eqref{pb_mu_def2}, with $i=1,...,k$,
is simple and is given by
$$\mu_i=\bigg((i-1)(i+N-2)\frac{4}{N(N-2)}+1\bigg)S^{2^*}_U,
$$
while
$$
\mu_{k+1}=\frac{S_U^{2^*}}{N(N-2)}\sqrt{(N-2)^2 + 4\lambda_1(D)}\bigg( 2 +\sqrt{(N-2)^2 + 4\lambda_1(D)}\bigg).
$$
Moreover, $\psi_i=R_i(r)Y_0(\theta)$ is an eigenfunction of \eqref{pb_mu_def2} corresponding to $\mu_i$, for any $i=1,...,k,$ where
$$
R_i=\frac{P_{i}(|x|^2)}{(1+|x|^2)^{\frac{N-2}{2}}}
$$
and $P_i(|x|^2)$ are given as in the Appendix \eqref{keigen}. Instead an eigenfunction corresponding to $\mu_{k+1}$ is given by  $\psi_{k+1}=R_{k+1}Y_1(\theta)$ where
$$R_{k+1}=(|x|^2+1)^{1-\frac{N}{2}-\beta(\lambda_1(D))}|x|^{\beta(\lambda_1(D))}$$
with
$2\beta(\lambda_1(D))=-(N-2) + \sqrt{(N-2)^2 + 4\lambda_1(D)} $ and $Y_1$ is an eigenfunction of \eqref{pb_lambda_def1} corresponding to $\lambda_1(D)$.
\end{proposition}

The proof is postponed after the following Corollary.

\begin{corollary}\label{corollary}
If $D$ is such that $\lambda_1(D)>N-1$ then
\begin{itemize}
\item[(i)] $\mu_1=S_U^{2^*}$ is simple with eigenfunction $U_{s}$; 
\item[(ii)] $\mu_2=(2^*-1)S_U^{2^*}$ is simple with eigenfunction $\partial_s U_{s}|_{s=1}$. 
\item[(iii)]$$\mu_3=\begin{cases}
\frac{S_U^{2^*}}{N(N-2)}\sqrt{(N-2)^2 + 4\lambda_1(D)}\bigg( 2 +\sqrt{(N-2)^2 + 4\lambda_1(D)}\bigg)&\text{for }N-1<\lambda_1(D)\le 2N\\
\\
\frac{(N+4)(N+2)}{N(N-2)}S_U^{2^*} \text{ (and is simple)}&\text{for }\lambda_1(D)> 2N.
\end{cases}
$$
Moreover if $N-1<\lambda_1(D)\le 2N$ a third eigenfunction $\psi_3$ is given by
$$\psi_3=(|x|^2+1)^{1-\frac{N}{2}-\beta(\lambda_1(D))}|x|^{\beta(\lambda_1(D))}Y_1
$$
where $2\beta(\lambda_1(D))=-(N-2) + \sqrt{(N-2)^2 + 4\lambda_1(D)} $ and $Y_1$ is an eigenfunction of \eqref{pb_lambda_def1} corresponding to $\lambda_1(D)$. Instead for $\lambda_1(D)> 2N$, the third eigenfunction is given by
$$\psi_3=(N-2(N+2)|x|^2+N|x|^4)(1+|x|^2)^{-\frac{N}{2}-1}.$$
\end{itemize}
\end{corollary}
\begin{proof}
    Note that if $\lambda_1(D)\in(N-1,2N)$ we have $k=2$ in Proposition \ref{lemma_autov}, while if $\lambda_1(D)>2N$ then \eqref{lemma_autov_cond} holds for some $k\ge 3$. Hence the thesis directly follows from Proposition \ref{lemma_autov}.
\end{proof}
\begin{proof}[Proof of Proposition \ref{lemma_autov}]
We study the eigenvalue problem \eqref{pb_mu_def2}.
We write the eigenvalue equation in \eqref{pb_mu_def2} by using spherical coordinates:
\beq\label{sphpb}
v''+(N-1)\frac{v'}{r}+\frac{1}{r^2}\Delta_{\S^{N-1}}v+\mu U^{2^*-2}v=0  \,,
\eeq
where $v'$ denotes the derivative of $v$ with respect to $r$ and $\Delta_{\S^{N-1}}$ denotes the Laplace-Beltrami operator on $\S^{N-1}$. 
From Proposition \ref{separation},  we can look for solutions of \eqref{sphpb} by using separation of variables : let $v=R(r)Y(\theta)$ be a solution to \eqref{sphpb}, then it holds
$$
Y(R''+(N-1)\dfrac{R'}{r})+\dfrac{1}{r^2}R\Delta_{\S^{N-1}}Y+\mu U^{2^*-2}RY=0\,,
$$
and hence, also taking into account that $\partial_\nu v= 0 $ on $\partial {\Sigma_D}$, we obtain that $Y$ is an eigenfunction of
\begin{equation}\label{angbp}
\begin{cases}
\Delta_{\S^{N-1}}Y+\lambda Y=0 & \text{on } D\\
\partial_{\nu_D} Y = 0 & \text{on } \partial D \,,
\end{cases}
\end{equation}
while $R$ satisfies
\beq\label{radpb}
\begin{cases}
R''+\dfrac{N-1}{r}R'- \dfrac{\lambda}{r^2}R+\mu U^{2^*-2}R=0 & \text{in }\R_+ \\
(r^{N-1}R'(r))(0)=0 & \\ 
R(r)=o(r^{2-N+\epsilon})  \quad \forall\epsilon> 0\text{ as }r\to + \infty \,. &
\end{cases}
\eeq

We consider \eqref{radpb}. To shorten the notation we write $\lambda_j$ instead of $\lambda_j(D)$. Thanks to the explicit expression of $U$, we can write the problem in the following Sturm Liouville form:
\begin{equation} \label{SL_eq}
\begin{cases}
(r^{N-1}R')'+r^{N-1}(-\lambda r^{-2}+\mu k_0^{2^*-2}(1+r^2)^{-2})R=0 & \text{in }\R_+ \\
(r^{N-1}R'(r))(0)=0 & \\ 
R(r)=o(r^{2-N+\epsilon})  \quad \forall\epsilon> 0\text{ as }r\to + \infty \,. &

\end{cases}
\end{equation}
Taking $\lambda=\lambda_j=\lambda_j(D)$, for $j=0,1,2,\ldots..$ we denote by $\mu_{k,\lambda_i}, k\in\N_+$, the sequence of eigenvalues of \eqref{SL_eq} obtained in correspondence of $\lambda_j$ and denote by $R^j_i$ the corresponding eigenfunction.

\textbf{Case $j=0$} i.e. $\lambda_0=0$.  The problem becomes:
 \begin{equation} \label{mu20}
\begin{cases}
(r^{N-1}R')'+r^{N-1}\mu_{k,\lambda_0} k_0^{2^*-2}(1+r^2)^{-2}R=0 & \text{in }\R_+ \\
(r^{N-1}R'(r))(0)=0 & \\ 
R(r)=o(r^{2-N+\epsilon})  \quad \forall\epsilon> 0\text{ as }r\to + \infty \,. &
\end{cases}
\end{equation}
From Appendix \ref{Appendix} we know that the eigenvalues are given by
\beq\label{lambdazerok}\mu_{k,\lambda_0}=\bigg((k-1)(k+N-2)\frac{4}{N(N-2)}+1\bigg)S^{2^*}\,,
\eeq
with corresponding eigenfunctions 
$$
R_{k,\lambda_0}=\frac{P_{k}(|x|^2)}{(1+|x|^2)^{\frac{N}{2}+k-2}}.
$$
In particular,  after some computations, we get that:
$$
R_{1,\lambda_0}=U_s, \quad R_{2,\lambda_0}=\partial_sU_s, \quad R_{3,\lambda_0}=\frac{N-2(N+2)|x|^2+N|x|^4}{(1+|x|^2)^{\frac{N}{2}+1}}.
$$

\textbf{Case $j=1$.} In this case we look for an eigenfunction of \eqref{SL_eq} of the following form
$$R_{\lambda_1(D)}=(1+r^2)^{1-\frac{N}{2}-\beta}r^\beta,
$$
for some $\beta\in \R$. From \eqref{SL_eq} we obtain that
$$
\bigg(r^{N-1}\bigg((1-\frac{N}{2}-\beta)(1+r^2)^{-\frac{N}{2}-\beta}2r^{\beta+1}+\beta(1+r^2)^{1-\frac{N}{2}-\beta}r^{\beta-1}\bigg)\bigg)'+
$$
$$
+r^{N-1}(-\lambda r^{-2}+\mu k_0^{2^*-2}(1+r^2)^{-2})(1+r^2)^{1-\frac{N}{2}-\beta}r^\beta=0
$$
that is
$$
r^{N-3+\beta}(1+r^2)^{-\frac{N}{2}-\beta-1}\bigg [(N-1)\bigg((1-\frac{N}{2}-\beta)(1+r^2)2r^{2}+\beta(1+r^2)^{2}\bigg)+
$$
$$
+\bigg((-2+N+2\beta)(N+2\beta)r^{4}+ (2\beta+1)(2-N-2\beta)(1+r^2)r^{2}+(\beta-1)\beta(1+r^2)^{2}\bigg)+
$$
$$
-\lambda_1(D) (1+r^2)^2+\mu k_0^{2^*-2}r^2\bigg]=0.
$$
Hence, $\beta$ must satisfy
$$
\begin{cases}
    r^4\bigg((2-N+\beta)(N-1-N-2\beta+2\beta+1)+\beta(N-1)+\beta(\beta-1)-\lambda_1(D)\bigg)=0\\
    r^2\bigg((N-1)(2-N)+(2\beta+1)(2-N-2\beta)+2(\beta-1)\beta-2\lambda+\mu k_0^{2^*-2}\bigg)=0\\
    \beta^2+\beta(N-2)-\lambda_1(D)=0,
\end{cases}
$$
i.e.,
$$
\begin{cases}
  N(2-N)-4(N-1)\beta-4\beta^2+\mu k_0^{2^*-2}=0\\
    \beta^2+\beta(N-2)-\lambda_1(D)=0,
\end{cases}
$$
and hence
$$
\begin{cases}
\mu k_0^{2^*-2}=(2\beta+N)(2\beta+N-2)\\
    \beta^2+\beta(N-2)-\lambda_1(D)=0.
\end{cases}
$$
This implies that
$$
R_{\lambda_1(D)}(r)=(r^2+1)^{1-\frac{N}{2}-\beta(\lambda_1(D))}r^{\beta(\lambda_1(D))} \,,
$$
with
\begin{equation}\label{eigencomp}
\begin{cases}
2\beta(\lambda_1(D))=-(N-2) + \sqrt{(N-2)^2 + 4\lambda_1(D)}  & \\
k_0^{\frac{4}{N-2}} \mu = \sqrt{(N-2)^2 + 4\lambda_1(D)}( 2 +\sqrt{(N-2)^2 + 4\lambda_1(D)}),  
\end{cases}
\end{equation}
and we notice that $R_{\lambda_1(D)}$ is an eigenfunction with no interior zeros, so it is the first one, i.e. it is $R_{1,\lambda_1}$ corresponding to $\mu_{1,\lambda_1}$. 
Since $k_0^{2^*-2}=S^{-2^*}N(N-2)$ (see Remark \ref{k0} below) we can write \eqref{eigencomp} in the following way:
\begin{equation}\label{betaeigenvalue}
\begin{cases}
2\beta(\lambda_1(D))=-(N-2) + \sqrt{(N-2)^2 + 4\lambda_1(D)}  & \\
\mu_{1,\lambda_1(D)}= \frac{S_U^{2^*}}{N(N-2)}\sqrt{(N-2)^2 + 4\lambda_1(D)}( 2 +\sqrt{(N-2)^2 + 4\lambda_1(D)}) .  
\end{cases}
\end{equation}
Note that $\mu_{1,\lambda_1(D)}$ is increasing in $\lambda_1(D)$.  

\textbf{Conclusion:} The conclusion follows by noticing that if $k\in\N_+$ is such that
$$
(k-1)(k+N-3)<\lambda_1(D)<k(k+N-2)
$$
then
$$
\mu_{k,\lambda_0}<\mu_{1,\lambda_1(D)}\le \mu_{k+1,\lambda_0}.
$$
From this computation, \eqref{lambdazerok} and \eqref{betaeigenvalue} we get the conclusion.
\end{proof}
\smallskip
\begin{remark}\label{k0}
We notice that $k_0^{2^*-2}=S_U^{-2^*}N(N-2)$, where $k_0$ is the constant such that $\|\nabla U\|_2=1$. Indeed from the definition of $U$ it holds
$$
\nabla U=-(N-2)k_0(1+|x|^2)^{-\frac{N}{2}}x.
$$
then
\beq\label{prop1k0}
\int_{\Sigma_D}(1+|x|^2)^{-N}|x|^2=\frac{1}{k_0^2 (N-2)^2}.
\eeq
Now \eqref{Sobquot} yields
\beq\label{prop2k0}
\int_{\Sigma_D} (1+|x|^2)^{-N}=\bigg(\frac{1}{S_U k_0}\bigg)^\frac{2N}{N-2}.
\eeq
Using radial coordinates and integrating by parts the radial terms we get
$$
\int_{\Sigma_D}(1+|x|^2)^{-N}|x|^2dx=\int_D d\omega\int_0^\infty \frac{r^{N+1}}{(1+r^2)^N}dr
$$
\beq \label{prop3k0}
=\frac{N}{2(N-1)}\int_D d\omega\int_0^\infty \frac{r^{N-1}}{(1+r^2)^{N-1}}dr=\frac{N}{2(N-1)}\int_{\Sigma_D} (1+|x|^2)^{-N+1}dx.
\eeq
Then from \eqref{prop1k0}, \eqref{prop2k0} and \eqref{prop3k0} we get
$$
\frac{1}{k_0^2 (N-2)^2}=\int_{\Sigma_D}(1+|x|^2)^{-N}|x|^2=\int_{\Sigma_D}(1+|x|^2)^{-N+1}-\int_{\Sigma_D}(1+|x|^2)^{-N}
$$
$$=\frac{2(N-1)}{N}\frac{1}{k_0^2 (N-2)^2}-(S_U k_0)^\frac{-2N}{N-2},
$$
namely
$$
k_0^{2^*-2}=S_U^{-2^*}N(N-2).
$$
\end{remark}

\section{Proof of the main results}\label{proofth1}
In this section we argue as in \cite{BE}. Let us start by observing that, for any $\varphi\in\mathcal{D}^{1,2}(\Sigma_D)$ and for any local minimizer $V$ with $\|\nabla V\|_2=1$ we have:
\beq\label{propdistanza}
d(\varphi,\mathcal{M_V})\le\|\nabla\varphi\|_2.
\eeq
This easily follows from the definition of the distance
$$
d(\varphi,\mathcal{M}_V)^2=\inf_{c,s}\bigg(\|\nabla \varphi\|_2^2+c^2-2c\int_{\Sigma_D} \nabla\varphi\cdot\nabla V_{s}dx \bigg)\,\le\inf_{s}\bigg(\|\nabla \varphi\|_2^2-2\bigg(\int_{\Sigma_D} \nabla\varphi\cdot\nabla V_{s}dx \bigg)^2\bigg)\,\le \|\nabla\varphi\|_2
$$
where we minimized w.r.t. $c$.

In order to prove Theorem \ref{thm_main}, we start with the following lemma.

\begin{lemma}\label{preliminarlemma}
Let $V$ be a nondegenerate local minimizer with $\|\nabla V\|_2=1$. Then any $\varphi\in\mathcal{D}^{1,2}(\Sigma_D)$ such that $d(\varphi,\mathcal{M}_V)<\|\varphi\|_2$ can be written as
$$
\varphi=c_0V_{s_{0}}+dv
$$
where $c_0,s_0$ are such that $$d:=d(\varphi,\mathcal{M}_V)^2=\|\nabla(\varphi-c_0V_{s_0})\|_2^2, \qquad c_0^2=\|\nabla\varphi\|_2^2-d^2,$$
and $v$ satisfies 
$$\|\nabla v\|_2=1,\qquad v\perp span[V_{s_0},\partial_sV_{s_0}], \qquad
\int_{\Sigma_D} V^{2^*-2}_{s_0}v^2dx\le\frac{1}{\mu_3}
$$
where $\mu_3$ is defined as in \eqref{pb_mu_def2}.
\end{lemma}
\begin{proof}
    We denote $\mathcal{M}_V$ simply by $\mathcal M$ and recall that it is a 2-dimensional manifold embedded in $\mathcal{D}^{1,2}(\Sigma_D)$:
$$
\R\setminus\{0\}\times\R^+\ni(c,s)\rightarrow cV_s\in \mathcal{D}^{1,2}(\Sigma_D).
$$
For any $\varphi\in\mathcal{D}^{1,2}(\Sigma_D)$ the distance of $\varphi$ to the manifold $\mathcal M$ is given by 
\begin{equation} \label{dist_def}
d(\varphi,\mathcal{M})^2=\inf_{c,s}\|\nabla(\varphi-cV_{s})\|_2^2=\inf_{c,s}\bigg(\|\nabla \varphi\|_2^2+c^2-2\int_{\Sigma_D} \nabla\varphi\cdot\nabla cV_{s}dx \bigg)\,,
\end{equation}
where we used the fact that $\|\nabla V_s\|_2 =1 $ for any $s \in R^+$. Since $d(\varphi,\mathcal{M})<\|\nabla \varphi\|_2$ then the infimum above is attained at a point $(c_0,s_0)\in\R\times\R^+$, with $c_0\ne 0$.  
Since $\mathcal M\setminus\{0\}$ is a smooth manifold we must have $(\varphi-c_0V_{s_0}) \perp T \mathcal M_{c_0V_{s_0}}$, which is the two-dimensional tangent space to $\mathcal{M}$ at the point $c_0V_{s_0}\in\mathcal{M}_V$. It is spanned by
\beq\label{deftg}
T\mathcal M_{c_0V_{s_0}}=span [V_{s_0},\partial_s V_{s}|_{s_0}].
\eeq
From Lemma \ref{lemma_spettro} we know that the spectrum of the operator $\mathcal{L}_{V_{s_0}}$ is discrete and by the min-max characterization we have
$$
\mu_3\le \frac{\int_{\Sigma_D} |\nabla \omega|^2dx}{\int_{\Sigma_D} V^{2^*-2}_{s_0}\omega^2dx}, \qquad \forall\omega\perp span [V_{s_0},\partial_s V_{s}|_{s_0}],
$$
since $V_{s_0}$ and $\partial_s V_{s}|_{s_0}$ are the first and second eigenfunction of the operator $\mathcal{L}_{V_{s_0}}$(see Proposition \ref{proposspettro}),
with equality attained if $\omega$ is the third eigenfunction. Then, from \eqref{deftg} we get 
\beq\label{mu3ineq}
\mu_3\le \frac{\int_{\Sigma_D} |\nabla \omega|^2dx}{\int_{\Sigma_D} V^{2^*-2}_{s_0}\omega^2dx}, \qquad \forall\omega\perp T\mathcal M_{c_0V_{s_0}}.
\eeq
From Proposition \ref{proposspettro} we know that $$\mu_1=S_{V_{s_0}}^{2^*}, \quad \mu_2=(2^*-1)S_{V_{s_0}}^{2^*},$$ and that all eigenvalues are independent of $s$.

Let $d=d(\varphi, \mathcal M)$. Since $(\varphi-c_0V_{s_0})\perp T\mathcal{M}_{c_0V_{s_0}}$ we can write 
$$
\varphi=c_0V_{s_0}+d v \,,
$$
where $v$ has norm $1$ in $\mathcal{D}^{1,2}(\Sigma_D)$ and is orthogonal to the tangent space $T\mathcal{M}_{c_0V_{s_0}}$. Hence, 
since $c_0V_{s_0}$ minimizes the distance of $\varphi$ from $\mathcal M$, we have the following orthogonality conditions
\begin{equation}
\int_{\Sigma_D} \nabla v \cdot \nabla V_{s_0} = 0 \quad \textmd{ and } \quad \int_{\Sigma_D} \nabla v \cdot \nabla ( \partial_s V_{s}|_{s_0} ) = 0,
\end{equation}
and from \eqref{dist_def} we obtain
\begin{equation} \label{dist2}
d^2 = \|\nabla \varphi\|_2^2+c_0^2-2\int_{\Sigma_D} \nabla\varphi\cdot\nabla c_0V_{s_0}dx = \|\nabla \varphi\|_2^2 -c_0^2 \,,
\end{equation}
where we have used that $\|\nabla V_{s_0}\|_2=1$. Moreover, since $V_{s_0}$ and $ \partial_s V_{s}|_{s_0} $ are eigenfunctions of $\mathcal{L}_{V_{s_0}}$, we also have
\begin{equation}\label{ort_cond}
\int_{\Sigma_D}  v \, V^{2^*-1}_{s_0} = 0 \quad \textmd{ and } \quad \int_{\Sigma_D} v \, V^{2^*-2}_{s_0} \, \partial_s V_{s}|_{s_0} = 0 \,.
\end{equation}

\end{proof}

\begin{proposition}\label{lemma1}
Let $V$ be a nondegenerate local minimizer, then the following inequality holds
$$
\|\nabla\varphi\|_2^2-S_V^{2}\|\varphi\|_{2^*}^2\ge d^2\bigg(1-\frac{\mu_2}{\mu_3}\bigg)+o(d^2),
$$
for all $\varphi\in\mathcal{D}^{1,2}(\Sigma_D)$ with $d(\varphi,\mathcal{M}_V)<\|\nabla \varphi\|_2.$ Moreover the coefficient of the second order term is sharp.
\end{proposition}
\begin{proof}
As before we assume that $\|\nabla V\|_2=1$ and by Lemma \ref{preliminarlemma} we can write $\varphi=c_0V_{s_0}+dv$.
An expansion in $d$ yields  
$$
|\varphi|^{2^*}=(c_0V_{s_0}+dv)^{2^*}=(c_0V_{s_0})^{2^*}+2^*(c_0V_{s_0})^{2^*-1}dv+\frac{2^*(2^*-1)}{2}(c_0V_{s_0})^{2^*-2}(dv)^2
+o(d^2),
$$
and then
\beq\label{expan}
\int_{\Sigma_D} |\varphi|^{2^*}dx=c_0^{2^*}S_V^{-2^*}+d2^*c_0^{2^*-1}\int_{\Sigma_D} V_{s_0}^{2^*-1}vdx+d^2\frac{2^*(2^*-1)}{2}c_0^{2^*-2}\int_{\Sigma_D} V^{2^*-2}_{s_0}v^2dx+o(d^2).
\eeq

From \eqref{ort_cond}, \eqref{mu3ineq} and $\|\nabla v\|_2=1$, we obtain
\beq\label{firstineq}
\int_{\Sigma_D} |\varphi|^{2^*}dx \le c_0^{2^*}S_V^{-2^*}+d^2c_0^{2^*-2}\frac{2^*(2^*-1)}{2}\frac{1}{\mu_3}+o(d^2)
\eeq

Thus, since $\mu_2=(2^*-1)S_V^{2^*}$, we get 
$$
\int_{\Sigma_D} |\varphi|^{2^*}dx \le c_0^{2^*}S_V^{-2^*}\bigg(1+d^2c_0^{-2}\frac{2^*}{2}\frac{\mu_2}{\mu_3}+o(d^2)\bigg) \,.
$$
and then
$$
\|\varphi\|_{2^*}^2\le c_0^2S_V^{-2}\bigg(1+d^2c_0^{-2}\frac{2^*}{2}\frac{\mu_2}{\mu_3}+o(d^2)\bigg)^\frac{2}{2^*}.
$$
Again, an expansion in $d$ yields
$$
\|\varphi\|_{2^*}^2\le c_0^2S_V^{-2}+d^2S_V^{-2}\frac{\mu_2}{\mu_3}+o(d^2),
$$
and from \eqref{dist2} we obtain  
\beq\label{cross}
\|\nabla\varphi\|_2^2-S_V^{2}\|\varphi\|_{2^*}^2\ge d^2\bigg(1-\frac{\mu_2}{\mu_3}\bigg)+o(d^2).
\eeq
Now we show that the coefficient of the second order term of \eqref{cross} is sharp. Let $\bar\varphi=V+\alpha w_3$ where $w_3$ is a third eigenfunction of $\mathcal{L}_V$ with $\|\nabla w_3\|=1$ and $\alpha>0$  is small. 
Now
\beq\label{3comp3}
d(\bar\varphi,\mathcal{M}_V)^2=\inf_{c,s}\bigg(\|\nabla \bar\varphi\|_2^2+c^2-2c\int_{\Sigma_D} \nabla\bar\varphi\cdot\nabla V_{s}dx \bigg)\,=\inf_{s}\bigg(\|\nabla\bar \varphi\|_2^2-\bigg(\int_{\Sigma_D} \nabla\bar\varphi\cdot\nabla V_{s}dx \bigg)^2\,\bigg),
\eeq
where the last equality follows from the minimization in $c$.
Since $w_3$ is a third eigenfunction and $V$ is a first one, then $w_3\perp V$and we have
$$
\|\nabla \bar\varphi\|_2^2= \|\nabla V\|_2^2+\alpha^2\|\nabla w_3\|_2^2=1+\alpha^2.
$$
Then from \eqref{3comp3} we get that 
$$
d(\bar\varphi,\mathcal{M})^2=\inf_{s}\bigg(1+\alpha^2-\bigg(\int_{\Sigma_D} \nabla (V+\alpha\omega_3)\cdot\nabla V_{s}dx \bigg)^2\,\bigg).
$$
From this, it is easy to see that for $\alpha$ small the closest point on $\mathcal M_V$ is $V$ and $d(\bar\varphi,\mathcal{M})=\alpha$.
Since $w_3$ is a third eigenfunction, then
$$
\mu_3= \frac{\int_{\Sigma_D} |\nabla w_3|^2dx}{\int_{\Sigma_D} V^{2^*-2}w_3^2dx},
$$
and from \eqref{expan} we have 
$$
\|\nabla\bar \varphi\|_2^2-S_V^{2}\|\bar \varphi\|_{2^*}^2 = d^2\bigg(1-\frac{\mu_2}{\mu_3}\bigg)+o(d^2),
$$
which implies that the constant of the second order term of \eqref{cross} is sharp and completes the proof.
\end{proof}
From Proposition \ref{lemma1} and Proposition \ref{lemma_autov} we obtain the following Corollary.
\begin{corollary}\label{corineq}
Let $\Sigma_D$ be a cone such that $\lambda_1(D)> N-1$, then the bubble $U$ is a nondegenerate local minimizer, and the following inequality holds
\beq\label{expanbubble}
\|\nabla\varphi\|_2^2-S_U^2\|\varphi\|_{2^*}^2\ge c_*d(\varphi,\mathcal{M})^2+o(d(\varphi,\mathcal{M})^2), 
\eeq
where
\beq\label{c^*}
c_*=\begin{cases}
    \frac{2\sqrt{N^2 + 4(\lambda_1(D)-N+1)}-2N + 4(\lambda_1(D)-N+1)}{\sqrt{N^2+ 4(\lambda_1(D)-N+1)}( 2 +\sqrt{N^2 + 4(\lambda_1(D)-N+1)})}& N-1<\lambda_1(D)\le 2N\\
\frac{4}{N+4}& \lambda_1(D)>2N
\end{cases}
\eeq
for all $\varphi\in\mathcal{D}^{1,2}(\Sigma_D)$ with $d(\varphi,\mathcal{M})<\|\nabla \varphi\|_2.$ 
\end{corollary}
\begin{proof}
From Proposition \ref{lemma1} we know that the thesis holds with 
\beq\label{c*comp}
c^*=1-\frac{\mu_2}{\mu_3}\eeq
where $\mu_2,\mu_3$ are respectively the second and the third eigenvalues of $\mathcal{L}_U$. From Proposition \ref{lemma_autov} we know that $\mu_2=(2^*-1)S_U^{2^*}$ and $$\mu_3=\begin{cases}
\frac{S^{2^*}}{N(N-2)}\sqrt{(N-2)^2 + 4\lambda_1(D)}\bigg( 2 +\sqrt{(N-2)^2 + 4\lambda_1(D)}\bigg)&\text{for }N-1<\lambda_1(D)\le 2N\\
\\
\frac{(N+4)(N+2)}{N(N-2)}S_U^{2^*} &\text{for }\lambda_1(D)> 2N.
\end{cases}
$$
substituting in \eqref{c*comp} we get the conclusion.
\end{proof}

\begin{remark}
Denoting $\mathcal{M}_U$ simply by $\mathcal{M}$, notice that if $\lambda_1(D)\le N-1$ then \eqref{expanbubble} is meaningless since $c_*$ becomes non positive. Indeed, if $\lambda_1(D)<N-1$ by arguing as in the proof of Lemma \ref{lemma1} and considering the spectral results obtained in Proposition \ref{lemma_autov}, we get the  following inequality
$$
\|\nabla\varphi\|_2^2-S_U^2\|\varphi\|_{2^*}^2\ge\bigg(1-\frac{(2^*-1)S_U^{2^*}}{\mu_2}\bigg)d(\varphi,\mathcal{M})^2+o(d(\varphi,\mathcal{M})^2), 
$$
where
$$
\mu_2=\frac{S_U^{2^*}}{N(N-2)}\sqrt{(N-2)^2 + 4\lambda_1(D)}( 2 +\sqrt{(N-2)^2 + 4\lambda_1(D)})\le(2^*-1)S_U^{2^*}
$$
for all $\varphi\in\mathcal{D}^{1,2}(\Sigma_D)$ with $d(\varphi,\mathcal{M})<\|\nabla \varphi\|_2$, and in this case we obtain that 
$$ \bigg(1-\frac{(2^*-1)S_U^{2^*}}{\mu_2}\bigg)<0.$$
If $\lambda_1(D)=N-1$ we have that $$\mu_2=(2^*-1)S_U^{2^*}$$ and then
$$
\|\nabla\varphi\|_2^2-S_U^2\|\varphi\|_{2^*}^2\ge o(d(\varphi,\mathcal{M})^2).
$$
As a matter of fact in \cite{CPP} it is proved that when $\lambda_1(D)<N-1$ the bubble cannot be a local minimizer.
\end{remark}

\begin{proof}[Proof of Theorem \ref{thm_main}]
\textbf{1)} We argue by contradiction. Suppose that there exists an infinite number of nondegenerate global minimizers $V^{(i)}, i \in\N$ with $\|\nabla V\|_2=1$. From \cite{CP} we know that up to a subsequence, it converges in $\mathcal{D}^{1,2}(\Sigma_D)$ to a minimizer $\bar V$ which must be nondegenerate and hence isolated. Then

$$
d(V^{(i)},\mathcal{ M}_{\bar V})\to 0 
$$
which contradicts that $\bar V $ is isolated.

\textbf{2)} We argue by contradiction. Assume that \eqref{quant_sobolev0} does not hold.  Then, using the Sobolev inequality 
$$
\|\nabla\varphi\|_2^2\ge S_D^2\|\varphi\|_{2^*}^2\quad \forall \varphi\in\mathcal{D}^{1,2}(\Sigma_D)
$$
we can find a sequence $\{\varphi_k\}, k\in\N$ such that
\begin{equation}\label{contradiction2}
\frac{\|\nabla\varphi_k\|_2^2- S_D^2\|\varphi_k\|_{2^*}^2}{\min_{i=1,...,m} d(\varphi_k,\mathcal{M}_{V^{(i)}})^2}\to 0,\quad  k\to\infty.
\end{equation}
By homogeneity, we can assume that $\|\nabla\varphi_k\|_2=1$. Since for all $i$ we have $d(\varphi_k,\mathcal{M}_{V^{(i)}})\le\|\nabla\varphi_k\|_2=1$, up to a subsequence, we can assume that $\min_i d(\varphi_k,\mathcal{M}_{V^{(i)}})\to \bar d\in[0,1]$. If $\bar d=0$ then Lemma \ref{lemma1} holds for some $\mathcal{M}_{V^{(i)}}$, and we get a contradiction. 
The other possibility is $0<\bar d\le 1$.  In this case, we must have
\begin{equation} \label{limit}
\|\nabla\varphi_k\|_2^2-\hat S^2\|\varphi_k\|_{2^*}^2\to 0,  
\end{equation}
which implies that $\varphi_k$ is a minimizing sequence, and then it converges to a minimizer $\bar V$ (see \cite{CP}), which belongs to $\mathcal{M}_{V^{(j)}}$ for some $j$. This implies 
$$
\min_id( \varphi_k,\mathcal{M}_{V^{(i)}})\to0 \quad\text{for }k\to \infty,
$$
which is a contradiction.
\end{proof}
\begin{corollary} 
Let ${\Sigma_D}$ be a cone as in \eqref{Sigma_def} such that the Sobolev quotient $Q_D$ admits a unique nondegenerate global minimizer $V$ up to rescaling and multiplication by a constant. Then the following  quantitative Sobolev inequality holds
\beq\label{quant_sobolev1}
\|\nabla\varphi\|_2^2- S_D^2\|\varphi\|_{2^*}^2\ge C_D d(\varphi,\mathcal{M}_{V})^2, \quad \forall \varphi\in\mathcal{D}^{1,2}(\Sigma_D)
\eeq
for some positive constant $C_D>0$. 
\end{corollary}
\begin{corollary}\label{corbub}
Let ${\Sigma_D}$ be a cone as in \eqref{Sigma_def} such that the only global minimizer, up to rescaling and multiplication by a constant, is the bubble $U$ and assume that $\lambda_1(D)> N-1$. Then the following  quantitative Sobolev inequality holds for all $\varphi\in\mathcal{D}^{1,2}(\Sigma_D)$
\beq
\|\nabla\varphi\|_2^2- S_D^2\|\varphi\|_{2^*}^2\ge C_D d(\varphi,\mathcal{M})^2, 
\eeq
where $0<C_D\le c^*$, where $c^*$ is given by \eqref{c^*}. 
\end{corollary}
As pointed out in the introduction, Corollary \ref{corbub} holds when the cone is convex or almost convex.
\begin{proof}
    The proof follows from Corollary \ref{corineq} and Theorem \ref{thm_main}.\end{proof}
By Proposition \ref{lemma1} we have that the sharp constant $C_D$ for the quantitative inequality defined in \eqref{bestconstdefic} is not greater than $\min_{i=1,...,m}\big(1-\frac{\mu_2^i}{\mu_3^i}\big)$. We conclude by showing that the inequality is strict.
\begin{theorem}\label{strict}
Let ${\Sigma_D}$ be a cone, with $\bar D\subset\S^{N-1}_+$, and assume that every minimizer is nondegenerate.
Let $V^{(1)},...,V^{(m)}$ and $C_D$ be as in \eqref{bestconstdefic}, then if either:
\begin{enumerate}
    \item $N=3,4,5$,
    \end{enumerate}
    or
    \begin{enumerate}
    \item[(2)] $N\ge 6$, $\Sigma_D$ is such that the only minimizers are the bubbles, $\lambda_1(D)\ge 2N$ and $D$ is symmetric w.r.t. two coordinates $x_i,x_j$ for $i,j\in\{1,...,N-1\}$, 
\end{enumerate}
then 
$$
C_D<\min_{i=2,...,m}\bigg(1-\frac{\mu_2^i}{\mu_3^i}\bigg).
$$
\end{theorem}
\begin{proof}
   Let $V$ be one of the minimizers $V^{(1)},...,V^{(m)}$, with $\|\nabla V\|_2=1$. We consider test functions of the form
$$
f_\varepsilon = V + \varepsilon \rho,
$$
for $\varepsilon$ small and $\rho \in T_{\mathcal{M}_V}^\perp$ satisfies a specific choice described below.

The orthogonality relation, and the fact that $-\Delta V = S_V^{2^*} V^{2^* - 1}$ easily imply
$$
\|\nabla f_\varepsilon\|_2^2 = \|\nabla V\|_{2}^2 + \varepsilon^2 \|\nabla \rho\|_{2}^2.
$$
A Taylor expansion in $\epsilon$ up to the fourth order yields
$$
(V + \varepsilon \rho)^{2^*} = V^{2^*} + \varepsilon 2^* V^{2^* - 1} \rho + \varepsilon^2 \frac{2^*(2^* - 1)}{2} V^{2^* - 2} \rho^2 +\varepsilon^3 \frac{2^*(2^* - 1)(2^*-2)}{6} V^{2^* - 3} \rho^3+
$$
$$+\varepsilon^4 \frac{2^*(2^* - 1)(2^*-2)(2^*-3)}{24} V^{2^* - 4} \rho^4+ o(\varepsilon^4),
$$
and using the orthogonality property $\int_{\mathbb{R}^d} V^{2^* - 1} \rho \, dx = 0$, we deduce
$$
\|V + \varepsilon \rho\|_{2^*}^{2} = \|V\|_{2^*}^{2} + \varepsilon^2 (2^* - 1) \|V\|_{2^*}^{2 - 2^*} \int_{\mathbb{R}^d} V^{2^* - 2} \rho^2 \, dx +\varepsilon^3 \frac{(2^* - 1)(2^* - 2)}{3} \|V\|_{2^*}^{2 - 2^*} \int_{\Sigma_D}V^{2^* - 3} \rho^3 \, dx +
$$

$$
+\varepsilon^4 \frac{(2^* - 1)(2^* - 2)(2^*-3)}{12} \|V\|_{2^*}^{2 - 2^*} \int_{\Sigma_D}V^{2^* - 4} \rho^4 \, dx+o(\varepsilon^4).
$$

Now taking $\rho$ as a third eigenfunction of the operator $\mathcal L_V$ with $\|\nabla\rho\|_2=1$ and remembering that $\|V\|_{2^*}=S_V^{-1}$, we have that the numerator of the quotient in \eqref{bestconstdefic} becomes
\begin{equation}\label{computationfinal}
\|\nabla f_\varepsilon\|_2^2 - S_D^2\|f_\varepsilon\|_{2^*}^2 = \bigg(1-\frac{\mu_2}{\mu_3}\bigg) \varepsilon^2 \|\nabla \rho\|_{2}^2 - \varepsilon^3 S_D^{2^*} \frac{(2^* - 1)(2^* - 2)  }{3}\int_{\Sigma_D} V^{2^* - 3} \rho^3 \, dx \end{equation}
$$- \varepsilon^4 S_D^{2^*} \frac{(2^* - 1)(2^* - 2) (2^*-3) }{12}\int_{\Sigma_D} V^{2^* - 4} \rho^4 \, dx + o(\varepsilon^4),
$$
where $\mu_2$ and $\mu_3$ are the second and the third eigenvalues of $\mathcal{L}_V$. As in the proof of Proposition \ref{lemma1} we have that for $\epsilon>0$ small enough, the minimum distance $\min_{i=1,...,m}\mathrm{dist}(f_\varepsilon, \mathcal{M}_{V^{(i)}})$ is achieved at $V$, hence
$$
\mathrm{dist}(f_\varepsilon, \mathcal{M}_V)^2 = \varepsilon^2 \|\nabla \rho\|_{2}^2= \varepsilon^2 .
$$

Hence, for $\varepsilon > 0$ small enough, by \eqref{computationfinal} we have that 
$$
\frac{\|\nabla f_\varepsilon\|_2^2-S_D^2\|f_\varepsilon\|_{2^*}^2}{\varepsilon^2} < 1-\frac{\mu_2}{\mu_3},
$$
provided that for $\varepsilon$ small enough
$$- \varepsilon^3 S_D^{2^*} \frac{(2^* - 1)(2^* - 2)  }{3}\int_{\Sigma_D} V^{2^* - 3} \rho^3 \, dx - \varepsilon^4 S_D^{2^*} \frac{(2^* - 1)(2^* - 2) (2^*-3) }{12}\int_{\Sigma_D} V^{2^* - 4} \rho^4 \, dx + o(\varepsilon^4)<0,
$$
that is
\begin{equation}\label{condfinal}
\int_{\Sigma_D} V^{2^* - 3} \rho^3 \, dx + \varepsilon  \frac{ (2^*-3) }{4}\int_{\Sigma_D} V^{2^* - 4} \rho^4 \, dx + o(\varepsilon)>0.
\end{equation}

\textbf{Case 1):} N=3,4,5. If $\rho$ is such that 
$$
\int_{\Sigma_D} V^{2^* - 3} \rho^3 \, dx\ne 0
$$
then we get the conclusion for $\epsilon$ small taking either $\rho$ or $-\rho$, which is still a third eigenfunction.  If $\rho$ is such that
\beq\label{3eigencond}
\int_{\Sigma_D} V^{2^* - 3} \rho^3 \, dx=0
\eeq
then we need
$$
\varepsilon  \frac{ (2^*-3) }{4}\int_{\Sigma_D} V^{2^* - 4} \rho^4 \, dx + o(\varepsilon)>0.
$$
This is true for $\epsilon$ small enough, because for $N=3,4,5$ we have $2^*>3$.

\textbf{Case 2):} $\lambda_1(D)\ge 2N$, the only minimizers are the bubbles and $D$ is symmetric with respect to two coordinates $x_i,x_j$, $i,j\in\{1,...,N-1\}$.
Then for $U$ as in \eqref{bubbledef} we know from the Appendix and Corollary \ref{corollary}, that the third eigenfunctions of $\mathcal{L}_U$ are given by $\rho(x)=U(x)\tilde\rho(\pi^{-1}(x))$ where $\pi$ is the stereographic projection and, $\tilde\rho$ is a homogeneous harmonic polynomial of degree two.
Then we call $\omega=\pi^{-1}(x)$ and we consider
$$
\tilde p(\omega)=\omega_i\omega_j+\omega_i\omega_l+\omega_j\omega_l
$$
where $\omega_l$ is any coordinate such that $l\ne i$ and $l\ne j$.

Then, by calling $\Omega=\pi^{-1}(\Sigma_D)$, we have that exists a constant $C(N)>0$ (depending only on $N$),
$$
\int_{\Sigma_D} V^{2^* - 3} \rho^3=C(N)\int_{\Omega} \tilde\rho(\omega)^3d\omega=C(N)\int_{\Omega} \omega_i^2\omega_j^2\omega_l^2>0
$$
since the monomial containing $\omega_i,\omega_i^3,\omega_j,\omega_j^3$ are odd and we suppose that the domain is symmetric with respect to those coordinates.
\end{proof}
\appendix
\section{An eigenvalue problem}\label{Appendix}
In this Appendix, we study the following eigenvalue problem:
 \begin{equation*} 
\begin{cases}
(r^{N-1}R')'+r^{N-1}\mu_{k,\lambda_0} k_0^{2^*-2}(1+r^2)^{-2}R=0 & \text{in }\R_+ \\
(r^{N-1}R'(r))(0)=0 & \\ 
R(r)=o(r^{2-N+\epsilon})  \quad \forall\epsilon> 0\text{ as }r\to + \infty \,. 
\end{cases}
\end{equation*}
We notice that studying this problem is the same as looking for the eigenvalues, corresponding to radial eigenfunctions, of the problem 
\beq\label{lambdazerobp}
\begin{cases}
-U^{2-2^*}\Delta v_k = \mu_{k,\lambda_0} v_k & \text{in } {\Sigma_D} \\
\partial_\nu v = 0 & \text{on } \partial {\Sigma_D} \\
v \in \mathcal{D}^{1,2}(\Sigma_D).
\end{cases}
\eeq
As in \cite[Section 2]{Premoselli},  we apply the transformation \beq\label{transfstereo}
\tilde v(y) =\frac{v}{U}(\pi(y))
\eeq
for $y\in \Omega\subset \S^N$, where $\pi$ is the stereographic projection from the North pole and $\Omega=\pi^{-1}({\Sigma_D})$.
In Cartesian coordinates, the projection from the north pole is given by:
$$
y_0=\frac{|x|^2-1}{|x|^2+1}\quad\text{ and }\quad y_i=\frac{2x_i}{|x|^2+1},
$$
where $x_1,...,x_N$ are the coordinates of $x\in{\Sigma_D}\subset\R^N$, and $y_0,...,y_N$ are the coordinates of  $y\in \Omega\subset \R^{N+1}$.
Since
$$
\bigg(\Delta_{g_0}+\frac{N(N-2)}{4}\bigg)\tilde v(y)=U(\pi(y))^{1-2^*}\Delta v(\pi(y))
$$
where $g_0$ is the metric on the sphere induced by the euclidean one in $\R^{N+1}$, then the equation \eqref{lambdazerobp} becomes
\beq\label{eqonsphere}
-\Delta_{g_0}\tilde v(y)=\frac{N(N-2)}{4}\bigg(S^{-2^*}\mu_{k,\lambda_0}-1\bigg)\tilde v, \qquad \text{in }\Omega\subset \S^N.
\eeq
This is the eigenvalue problem for the Laplacian on the set $\Omega$ on the sphere. We are looking for radial eigenfunction of \eqref{lambdazerobp} that through the transformation \eqref{transfstereo}, are the eigenfunctions of the Laplacian on $\Omega$ that depend only on the $y_0$ coordinate. Since we are looking for the radial one and we have the Neumann boundary condition, we are interested in the eigenfunction of $-\Delta_{g_0}$ that depends on $y_0$ on the whole sphere $\S^N$. It is known that if we have a harmonic, homogeneous, polynomial of degree $k-1$ for $k\in\N_+$, then it is an eigenfunction corresponding to the eigenvalue $(k-1)(k+N-2)$.  From \cite[pag. 85]{ABR} we also know that the $(k-1)$-th partial derivative of $|y|^{1-N}$ will give us an harmonic homogeneous polynomial of degree $k-1$. Then if we consider
$$
\partial^{k-1}_{y_0}|y|^{1-N}
$$
we obtain a radial eigenfunction corresponding to the $(k-1)$-th eigenvalue.

For $k=1$ we have
$$
\tilde V^{(1)}(x)=|y|^{1-N}=1,
$$
for $k=2$
$$
\tilde v_2(x)=\partial_{y_0}|y|^{1-N}=(1-N)|y|^{-N-1}y_0=2(1-N)\frac{|x|^2-1}{|x|^2+1},
$$
for $k=3$
$$
\tilde v_3(x)=\partial^{2}_{y_0}|y|^{1-N}=(N-1)(N+1)|y|^{-N-3}y_0^2-(N-1)|y|^{-N-1}=
$$
$$
=(N-1)\frac{N(|x|^4+1)-2(N+2)|x|^2}{(|x|^2+1)^2}.
$$
Then from the transformation $\tilde v(y) =\frac{v}{U}(\pi(y))$ we get that, up to a constant,
$$
V^{(1)}=U, \quad v_2=k_0(1+|x|^2)^\frac{-N}{2}(1-|x|^2)=\partial_sV_{s=1}, 
$$
$$
v_3=k_0(1+|x|^2)^\frac{-N-2}{2}(N|x|^4-2(N+2)|x|^2+N).
$$
and more generally
\beq\label{keigen}
v_k U^{-1}= (\partial^{k}_{y_0}|y|^{1-N})(\pi^{-1}(x)):=P_i(|x|^2).
\eeq
Then we get that the eigenvalues $\mu_{k,\lambda_0}$ of \eqref{lambdazerobp} satisfy
$$
\frac{N(N-2)}{4}\bigg(S^{-2^*}\mu_{k,\lambda_0}-1\bigg)=(k-1)(k+N-2)
$$
and then they are given by
\beq\mu_{k,\lambda_0}=\bigg((k-1)(k+N-2)\frac{4}{N(N-2)}+1\bigg)S^{2^*}.
\eeq
In particular, for $k=1$ we obtain
$$\mu_{1,\lambda_0}=S^{2^*},
$$
for $k=2$
$$\mu_{2,\lambda_0}=(2^*-1)S^{2^*},
$$
and for $k=3$
$$\mu_{3,\lambda_0}=\frac{(N+4)(N+2)}{N(N-2)}S^{2^*}.
$$

\end{document}